\documentclass[12pt]{article}
\usepackage{multirow}
\usepackage[english]{babel}
\usepackage[utf8]{inputenc}
\usepackage[OT1]{fontenc}
\usepackage{amsfonts, amsmath, amsthm, amssymb}
\usepackage{graphicx}
\usepackage{setspace}
\usepackage{listings}
\usepackage[margin=1in]{geometry}
\usepackage{xcolor}
\newcounter{countCode}
\lstnewenvironment{code} [1][caption=Ponme caption, label=default]{%
	 
	\setcounter{lstlisting}{\value{countCode}} 
	\lstset{ %
	language=matlab,
	basicstyle=\ttfamily\footnotesize,      
	numbers=left,                  
	numberstyle=\sc,      
	stepnumber=1,                  
	numbersep=5pt,                
	numberstyle=\color{red!50!blue},
    	backgroundcolor=\color{lightgray!20},
	rulecolor=\color{blue},
	keywordstyle=\color{red}\bfseries,
	showspaces=false,               
	showstringspaces=false,         
	showtabs=false,                
	frame=single,                   
	framexleftmargin=0mm,
	numberblanklines=false,
	xleftmargin=5pt,
	breaklines=true,
	breakatwhitespace=true,
	breakautoindent=true,
	captionpos=t,
	texcl=true,
	tabsize=2,                      
	extendedchars=true,
	inputencoding=utf8, 
	escapechar=\%,
	morekeywords={print, println, size, background, strokeWeight, fill, line, rect, ellipse, triangle, arc, save, PI, HALF_PI, QUARTER_PI, TAU, TWO_PI, width, height,},
	emph=[1]{print,println,}, emphstyle=[1]{\color{blue}}, 
	emph=[2]{width,height,}, emphstyle=[2]{\bf\color{violet}}, 
	emph=[3]{PI, HALF_PI, QUARTER_PI, TAU, TWO_PI}, emphstyle=[3]\color{orange!50!violet}, 
	emph=[4]{line, rect, ellipse, triangle, arc,}, emphstyle=[4]\color{green!70!black}, 
	#1}}{\addtocounter{countCode}{1}}

\title{The Set-Theoretic Form of Twin Prime Distribution and Its Odd-Even Imbalance}
\author{HaoJie Huang}

\date{
    \begin{center}
        School of Energy and Power Engineering, University of Shanghai for Science and Technology, Shanghai 200093, PR China \\
        Email: hjhuang@usst.edu.cn
    \end{center}
}

\begin{document}
\begin{spacing}{1.5}%%行间距变为double-space
\begin{spacing}{1}%%行间距变为double-space
\maketitle \tableofcontents

\begin{abstract}
The prime number problem falls within the realm of number theory, specifically elementary number theory. Current research approaches have unnecessarily complicated this matter. In contrast to more advanced mathematical tools, the methods of elementary number theory can effectively address the twin prime problem. The primary contribution of this article lies in establishing a set that systematically includes all twin prime pairs without omissions. This set's distribution is governed by a precise general solution formula. Furthermore, an analysis of the distribution set reveals characteristics of parity balance, enabling us to determine whether twin prime pairs are finite. Finally, our method can be extended to the general case of twin primes, such as the Polignac's conjecture.
\end{abstract}
\end{spacing}

\section {Rediscovering Prime Numbers}
\label{sec:Rediscovering Prime Numbers}

 Twin Prime Conjecture: There exist infinitely many prime numbers $p$ such that $p+2$ is also a prime number. For example, (3,5), (5,7), (11,13), (17,19), and so on. Before we formally embark on this conjecture, let's first understand what prime numbers are.

\subsection {Normal definition of prime numbers} 
The simple definition of a prime number is that it is not divisible by any integer except number 1 and itself. And, we see that of whole the prime numbers, all are odd except the number 2. 

 \subsection {Prime numbers defined in terms of sets}
If any odd number $(2z+1)$ can be expressed in the following form:
\begin{equation}
   {(2z+1)}={(2x+1)}\cdot{(2y+1)},{x}\in {1 \sim k_{1}},{y}\in {1 \sim k_{2}}.
   \label{1}
\end{equation}
$k_{1}$ and $k_{2}$ are positive integers. Then $(2z+1)$ is not a prime number.We expand equation~(\ref{1}) above:
\begin{equation}
   {z}={2xy+x+y},{x}\in {1 \sim k_{1}},{y}\in {1 \sim k_{2}}.
   \label{2}
\end{equation}
Equation~(\ref{2}) gives us a simple set of all numbers that are not prime numbers. For example, when $x=1$ and $y=1$, thus $z=4$, so number $9$ is not prime. For example, when $x = 2$ and $y = 1$, thus $z = 7$, so number $15$ is not prime. And any other value that cannot be expressed by $2xy+x+y$, such as $z$ = 1, 2, 3, 5, 6.... , this $(2z+1)$ is a prime number. In particular, we are discussing this value of $z$ in all subsequent expressions, not $(2z+1)$.

\section {Rediscovering Twin Primes}
\label{sec:Rediscovering Twin Primes}

And the so-called twin prime numbers, i.e., it is sufficient to satisfy that both $(2z+1)$ and $(2(z+1)+1)$ are prime. We note that in equation~(\ref{2}) obtained earlier, if we bring $y = 1$ to $2xy+x+y$, we get this interesting sequence: $[3x + 1]$. If you remember the definition of a twin prime, you will understand that all cases of twin primes will occur between neighbouring elements in the sequence $[3x+1]$. This simply means that in the sequence $[4,7,10,13,16,19,22...]$, all twin prime numbers occur between two neighbouring elements. For example: 5 and 6, 8 and 9, 11 and 12... and so on. 

Of course, not all twin prime numbers are satisfied between every two neighbouring elements of the $[3x+1]$ sequence due to the effect of other more widely spaced sequence values which generated when $y$ takes on larger values in equation~(\ref{2}). For example, 11 and 12 (12 is in the $[5x+2]$ sequence) represent 23 and 25, which are not twin prime pairs. The above information can be better understood from the table 1. In fact, one of the more interesting things about this $[3x+1]$ sequence is that it starts at 4. The number 1, 2, and 3 that precede it represent 3, 5, and 7, which are obviously also twin prime numbers, and the only set of triple-twin prime numbers.

\begin{table}
\begin{center}
\resizebox{\textwidth}{!}{%
\scriptsize
\begin{tabular}{lccc}
\hline
Sequence {[3x+1]}        & Natural number correspondence           & Property                   \\ \hline
{\color[HTML]{FE0000} 4} & {\color[HTML]{FE0000} 9} & {\color[HTML]{FE0000} Composite number} \\
5                        & 11                        &                            \\
6                        & 13                        & \multirow{-2}{*}{Twin prime pair}     \\
{\color[HTML]{FE0000} 7} & {\color[HTML]{FE0000} 15} & {\color[HTML]{FE0000} Composite number} \\
8                        & 17                        &                            \\
9                        & 19                        & \multirow{-2}{*}{Twin prime pair}     \\
{\color[HTML]{FE0000} 10} & {\color[HTML]{FE0000} 21} & {\color[HTML]{FE0000} Composite number} \\
11                        & 23                        &                            \\
12                        & 25                        & \multirow{-2}{*}{No Twin prime pair}     \\
{\color[HTML]{FE0000} 13} & {\color[HTML]{FE0000} 27} & {\color[HTML]{FE0000} Composite number}\\
\hline
\end{tabular}
}
\end{center}
\caption{Elements and neighbouring elements in $3x+1$ sequences corresponding to twin prime pairs
}
\end{table}

\section{General Solution Formula}

Our current task is to identify which elements in the $[3x+1]$ sequence belong to twin primes and which do not, and to discern the underlying patterns. We are revising our conjecture about twin primes into a more specific mathematical expression:

\begin{equation}
\begin{cases}
{3x_{1}+1=p_{1}},\text{positive integer solutions},[3x+1] \\

\left. \begin{aligned}
{(2n_{1}+1)x_{2}+n_{1}=p_{2}=p_{1}+1},\text{no positive integer solutions}\\
{(2n_{2}+1)x_{3}+n_{2}=p_{3}=p_{1}+2},\text{no positive integer solutions}\\
\end{aligned} \right\} \Rightarrow \text{twin primes!}
\\
{3x_{1}+4=p_{4}=p_{1}+3},\text{positive integer solutions},[3x+4]
    \end{cases}
    \label{3}
\end{equation}

Thus,
\begin{equation}
\begin{cases}
{x_{2}=\frac{3x_{1}+2-n_{1}}{2n_{1}+1}},\\
{x_{3}=\frac{3x_{1}+3-n_{2}}{2n_{2}+1}}.
    \end{cases}
    \label{4}
\end{equation}

As long as either of the expressions for  $x_{2}$ and $x_{3}$ has a general solution, the requirement for twin prime pairs is not satisfied. For the sake of convenience in generalization, we uniformly rewrite $n_{2}$ and $n_{3}$ as $n$, which does not affect the final outcome. 

For example:

When $x_{2}$ equals 1 and $x_{3}$ equals 1, $x_{1}$ in the equation~(\ref{4}) can be obtained as follows,
\begin{equation}
\begin{cases}
{x_{1}={n}-\frac{1}{3}},\text{no general solution}\\
{x_{1}={n}-\frac{2}{3}},\text{no general solution}.
    \end{cases}
    \label{5}
\end{equation}

When $x_{2}$ equals 2 and $x_{3}$ equals 2,
\begin{equation}
\begin{cases}
{x_{1}={5n}},\\
{x_{1}={5n-2}}.
    \end{cases}
    \label{6}
\end{equation}

When $x_{2}$ equals 3 and $x_{3}$ equals 3,
\begin{equation}
\begin{cases}
{x_{1}={7n-2}},\\
{x_{1}={7n}}.
    \end{cases}
    \label{7}
\end{equation}

When $x_{2}$ equals 4 and $x_{3}$ equals 4, there is no general solution.

Similarly, we can continue to obtain the subsequent general solution expressions, thereby organizing and obtaining:
\begin{equation}
\begin{cases}
{x_{1}={5n-2};{11n-3};{17n-4};{23n-5};...},\\
{x_{1}={5n+0};{11n+1};{17n+2};{23n+3};...},\\
{x_{1}={7n-2};{13n-3};{19n-4};{25n-5};...},\\
{x_{1}={7n+0};{13n+1};{19n+2};{25n+3};...}.
    \end{cases}
    \label{8}
\end{equation}
Summing up, we arrive at the final generalized expression for the solution about $x_{1}$:

\begin{equation}
  x_{1}=\begin{cases}
    (6m-1)n\pm{m}-1\\
    (6m+1)n\pm{m}-1,{m}\in {1 \sim k_{1}},{n}\in {1 \sim k_{2}}.
  \end{cases}
  \label{9}
\end{equation}
$k_{1}$ and $k_{2}$ are positive integers. 

The significance of this general solution formula lies in the fact that within the sequence $[3x+1]$, all numerical points satisfying the general solution formula $x_{1}$, and their subsequent two values ($3x_{1}+2$ and $3x_{1}+3$, requiring further transformation to odd numbers), are guaranteed not to be twin prime numbers. Conversely, if they do not satisfy the general solution formula, they are guaranteed to be twin prime numbers.

Here's a simple illustration. We will partially display the numerical points satisfying the general solution formula $x_{1}$ in the form of a table below. For example, if 1, 2, 4, 6, and 9 do not appear in this table, then the corresponding values for the $[3x+1]$ sequence are 4, 7, 13, 19, 28. Reverting these values back to their corresponding odd numbers, we can get the numbers: 9, 15, 27, 39, 57. The two odd numbers after these numbers must be twin prime numbers. For example: (11, 13), (17, 19), (29, 31), (41, 43), (59, 61). You can verify the situation of other twin prime pairs yourself.

\begin{table}[]
\centering
\label{tab:my-table}
\resizebox{0.2\columnwidth}{!}{%
\begin{tabular}{|llll|l|}
\hline
\multicolumn{1}{|l|}{\color[HTML]{FE0000} \textbf{3}}  & \multicolumn{1}{l|}{\textbf{8}}                         & \multicolumn{1}{l|}{\textbf{13}}                        & \textbf{18}                         &                                 \\ \cline{1-4}
\multicolumn{1}{|l|}{\color[HTML]{FE0000} \textbf{5}}  & \multicolumn{1}{l|}{\textbf{10}}                        & \multicolumn{1}{l|}{\textbf{15}}                        & \textbf{20}                         &                                 \\ \cline{1-4}
\multicolumn{1}{|l|}{\color[HTML]{FE0000} \textbf{5}}  & \multicolumn{1}{l|}{\textbf{12}}                        & \multicolumn{1}{l|}{\textbf{19}}                        & \textbf{26}                         &                                 \\ \cline{1-4}
\multicolumn{1}{|l|}{\color[HTML]{FE0000} \textbf{7}}  & \multicolumn{1}{l|}{\textbf{14}}                        & \multicolumn{1}{l|}{\textbf{21}}                        & \textbf{28}                         &                                 \\ \cline{1-4}
\multicolumn{1}{|l|}{\textbf{8}}  & \multicolumn{1}{l|}{{\color[HTML]{FE0000} \textbf{19}}} & \multicolumn{1}{l|}{\textbf{30}}                        & \textbf{41}                         &                                 \\ \cline{1-4}
\multicolumn{1}{|l|}{\textbf{10}} & \multicolumn{1}{l|}{{\color[HTML]{FE0000} \textbf{23}}} & \multicolumn{1}{l|}{\textbf{36}}                        & \textbf{49}                         &                                 \\ \cline{1-4}
\multicolumn{1}{|l|}{\textbf{12}} & \multicolumn{1}{l|}{{\color[HTML]{FE0000} \textbf{23}}} & \multicolumn{1}{l|}{\textbf{34}}                        & \textbf{45}                         &                                 \\ \cline{1-4}
\multicolumn{1}{|l|}{\textbf{14}} & \multicolumn{1}{l|}{{\color[HTML]{FE0000} \textbf{27}}} & \multicolumn{1}{l|}{\textbf{40}}                        & \textbf{53}                         &                                 \\ \cline{1-4}
\multicolumn{1}{|l|}{\textbf{13}} & \multicolumn{1}{l|}{\textbf{30}}                        & \multicolumn{1}{l|}{{\color[HTML]{FE0000} \textbf{47}}} & \textbf{64}                         &                                 \\ \cline{1-4}
\multicolumn{1}{|l|}{\textbf{15}} & \multicolumn{1}{l|}{\textbf{34}}                        & \multicolumn{1}{l|}{{\color[HTML]{FE0000} \textbf{53}}} & \textbf{72}                         &                                 \\ \cline{1-4}
\multicolumn{1}{|l|}{\textbf{19}} & \multicolumn{1}{l|}{\textbf{36}}                        & \multicolumn{1}{l|}{{\color[HTML]{FE0000} \textbf{53}}} & \textbf{70}                         &                                 \\ \cline{1-4}
\multicolumn{1}{|l|}{\textbf{21}} & \multicolumn{1}{l|}{\textbf{40}}                        & \multicolumn{1}{l|}{{\color[HTML]{FE0000} \textbf{59}}} & \textbf{78}                         &                                 \\ \cline{1-4}
\multicolumn{1}{|l|}{\textbf{18}} & \multicolumn{1}{l|}{\textbf{41}}                        & \multicolumn{1}{l|}{\textbf{64}}                        & {\color[HTML]{FE0000} \textbf{87}}  &                                 \\ \cline{1-4}
\multicolumn{1}{|l|}{\textbf{20}} & \multicolumn{1}{l|}{\textbf{45}}                        & \multicolumn{1}{l|}{\textbf{70}}                        & {\color[HTML]{FE0000} \textbf{95}}  &                                 \\ \cline{1-4}
\multicolumn{1}{|l|}{\textbf{26}} & \multicolumn{1}{l|}{\textbf{49}}                        & \multicolumn{1}{l|}{\textbf{72}}                        & {\color[HTML]{FE0000} \textbf{95}}  &                                 \\ \cline{1-4}
\multicolumn{1}{|l|}{\textbf{28}} & \multicolumn{1}{l|}{\textbf{53}}                        & \multicolumn{1}{l|}{\textbf{78}}                        & {\color[HTML]{FE0000} \textbf{103}} & \multirow{-16}{*}{\textbf{...}} \\ \hline
\multicolumn{4}{|l|}{\textbf{...}}                                                                                                                                                          & \textbf{...}                    \\ 
\hline
\end{tabular}%
}
\caption{Partial values of the array matrix $x_{1}$}
\end{table}

\section{Odd-even Distribution Properties}
The next step is to study the parity distribution characteristics of the numerical points in this set. We know that if a set of pairs of twin prime numbers is finite, assuming we have found the largest pair, then the values in the matrix of generic solution array corresponding to numbers after this maximum number will continuously cover the subsequent integer domain. At this point, the counts of odd and even numbers should be consistent. That is, as tends to infinity, the difference in the number of all odd and even numbers in the matrix of generic solution array will tend to a fixed value or what is called parity balance. Thus, analyzing the balance of odd and even numbers is crucial for determining whether a sequence of twin prime numbers is infinite or finite.

Before analysing it I would like you to notice some features of the entire matrix of the generic solution array:\\
(a) The arrays (groups of four elements) in the matrix of generic solution array are alternately parity, whereby we can assume that the distribution of the number of parities is homogeneous, taking into account repetitions.\\
(b) Apart from the symmetric axis of the entire array matrix (highlighted in red), it is difficult to find stable repeated array sequences; And the values of the symmetric axes are all odd numbers.\\
(c) Apart from the axis of symmetry, other repeated values include both odd and even numbers, with no clear advantage in odd or even parity.\\
(d) Dividing the array matrix by its axis of symmetry, where the upper and lower parts are completely symmetrical, we consider that half of the symmetric axis array should be duplicated, thus ensuring that half of the entire solution array matrix is duplicated.

Based on this, we believe that when the maximum value of the array matrix, denoted as $x_{max}=6n^2+2n-1$, there is a unique set of sequences that is stably repetitive and odd in quantity, totaling approximately $2n$. However, within this set, approximately $0.2n$ odd numbers do not repeat well (which can be easily demonstrated: the tail numbers are 5, 3, 3, 5, 9, where except for 9, the tails of 5, 3, 3, 5 repeat the tails of the first two rows of the array matrix at the very beginning). That is to say, there are fewer even numbers than odd numbers in the array matrix because odd numbers lack stable repetition along the symmetry axis, resulting in a deficit of $0.2n$. Therefore, we can infer that in the sequence of numbers outside the array matrix range (the number corresponding to the twin prime number pair), there are more even numbers than odd numbers.\\

\begin{figure}
    \centering
    \includegraphics[width=1\linewidth]{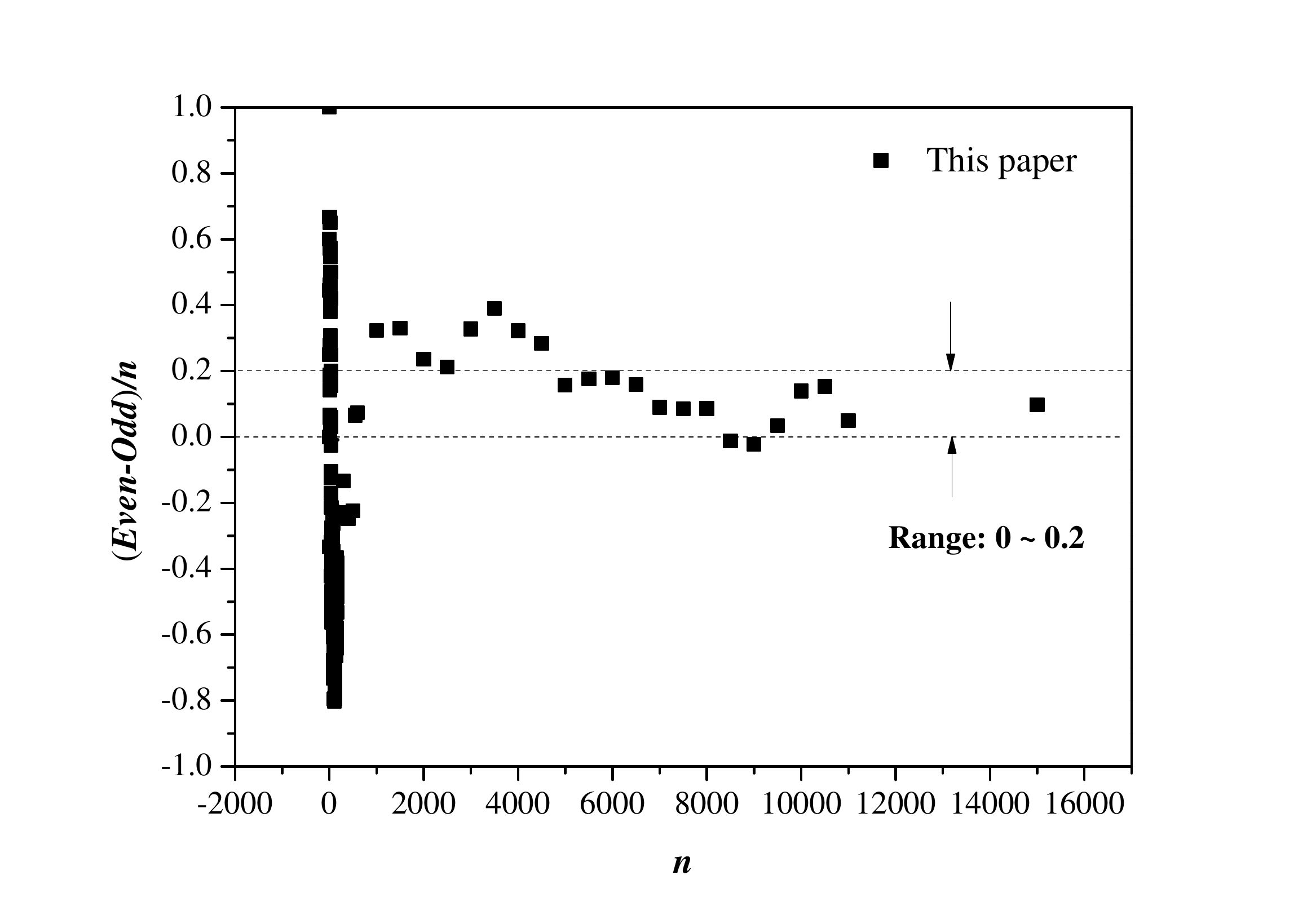}
    \caption{Relative differences between even and odd numbers obtained in different $n$ ranges}
\end{figure}

Here, we present the results of the analysis through Figure 1. The horizontal axis represents the values of $n$, while the vertical axis represents the ratio of the difference between the count of all even numbers and all odd numbers, which are not within the array matrix and are less than $x_{max}$, to $n$. It can be observed that as $n$ tends to a large value, the values on the vertical axis fluctuate within the range of 0 to 0.2. Based on our preliminary calculations, it aligns with our speculation ($\frac{Even-Odd}{n}>0$, where the number of even numbers is greater than the number of odd numbers), and the results fall within the range of 0 to 0.2, which is consistent. To summarise, the infinity of pairs of twin prime numbers can be illustrated by the characteristic of parity imbalance.\\

\noindent\textbf{Matlab Program} \\
\begin{spacing}{1}%%行间距变为double-space
Matlab program to generate twin prime pairs using array matrices:
\begin{code}[caption=Matlab code, label=default]
clc
clear all
n=150;(you can set a larger number)
k=n;
nn=6*n*n+2*n-1;
m=floor((nn+2)/5);
n=floor((nn+2)/5);
k1=0;
k2=0;
k3=0;
k4=0;
for i=1:1:m
    for j=1:1:n
        p1=(6*i+1)*j+i-1;
        if p1<=nn
            k1=k1+1;
            pp1(k1)=p1;
        end
        p2=(6*i-1)*j+i-1;
        if p2<=nn
            k2=k2+1;
            pp2(k2)=p2;
        end
        p3=(6*i+1)*j-i-1;
        if p3<=nn
            k3=k3+1;
            pp3(k3)=p3;
        end
        p4=(6*i-1)*j-i-1;
        if p4<=nn
            k4=k4+1;
            pp4(k4)=p4;
        end        
    end
end
pp=[pp1';pp2';pp3';pp4'];
p=sort(pp(:));
p=unique(p);
q=1:nn;
C = setdiff(q,p);
C_twin_prime1=(C'*3+1)*2+1+2;
C_twin_prime2=(C'*3+1)*2+1+4;
result = [C_twin_prime1, C_twin_prime2];
disp('Twin Prime:');
disp(result);
[ind]=find(rem(C,2)==0);
disp('(Even-Odd)/n=')
disp((-length(C)+length(C(ind))+length(C(ind)))/k)

\end{code}
\end{spacing}

\section{Expand to the Polignac's Conjecture}
Polignac's conjecture: There exist infinitely many prime numbers $p$, such that $p+2k$ is also prime, for any positive integer $k$.

When $k=1$, it degenerates into the situation of twin primes: $p$ and $p+2$.
The generalized expression for the solution regarding $x_{1}$ can be found in equation~(\ref{9}).

When $k=2$, namely $p$ and $p+4$, we can apply the same approach to obtain $x_{1}$.

\begin{equation}
\begin{cases}
{3x_{1}+1=p_{1}},\text{positive integer solutions},[3x+1] \\

\left. \begin{aligned}
{(2n_{1}+1)x_{2}+n_{1}=p_{2}=p_{1}+2},\text{no positive integer solutions}\\
{(2n_{2}+1)x_{3}+n_{2}=p_{3}=p_{1}+4},\text{no positive integer solutions}\\
\end{aligned} \right\} \Rightarrow \text{Twin primes of the form $2k$($k=2$)}
\\
{3x_{1}+4=p_{4}=p_{1}+3},\text{positive integer solutions},[3x+4]
    \end{cases}
    \label{10}
\end{equation}

\begin{equation}
  x_{1}=\begin{cases}
    (6m-1)n\pm\frac{1}{2}({2m-1})-\frac{3}{2}\\
    (6m+1)n\pm\frac{1}{2}({2m+1})-\frac{3}{2},{m}\in {1 \sim k_{1}},{n}\in {1 \sim k_{2}}.
  \end{cases}
  \label{11}
\end{equation}
$k_{1}$ and $k_{2}$ are positive integers. 

When $k=3$, two cases arise, and the solution method is similar. Without further elaboration, we directly provide the results:
\begin{equation}
  x_{1}=\begin{cases}
    (6m-1)n+\frac{1}{2}({2m-3})\pm\frac{1}{2}\\
    (6m+1)n-\frac{1}{2}({2m+3})\pm\frac{1}{2},{m}\in {1 \sim k_{1}},{n}\in {1 \sim k_{2}}.
  \end{cases}
  \label{11}
\end{equation}
$k_{1}$ and $k_{2}$ are positive integers. \\
Alternatively:
\begin{equation}
  x_{1}=\begin{cases}
    (6m-1)n-\frac{1}{2}({2m+3})\pm\frac{1}{2}\\
    (6m+1)n+\frac{1}{2}({2m-3})\pm\frac{1}{2},{m}\in {1 \sim k_{1}},{n}\in {1 \sim k_{2}}.
  \end{cases}
  \label{11}
\end{equation}
$k_{1}$ and $k_{2}$ are positive integers. 

For other values of $k$, readers can derive the solutions themselves. As for the existence of a comprehensive general solution formula applicable to all $k$, I believe it's possible, and I leave it to everyone to complete.

\section {Reference}
\label{sec:Reference}

This article does not require references, but thanks to individuals including Euclid, the ancient Greek mathematician.

\section {Acknowledgements}
\label{sec:Acknowledgements}
The article is heaven-born, the clumsy hand occasionally obtains it.\\
Hereby dedicate this article to my disappointed self.

\end{spacing}
\end{document}